%% file: anisotropic.tex
\documentclass[11pt,twoside]{article}
\usepackage{amsmath,amsfonts,amscd,amsthm}
\newtheorem{thm}{Theorem}[section]
\newtheorem{prop}[thm]{Proposition}
\newtheorem{lem}[thm]{Lemma}
\newtheorem{cor}[thm]{Corollary}

\theoremstyle{definition}
\newtheorem{defin}[thm]{Definition}

\theoremstyle{remark}
\newtheorem{rem}[thm]{Remark}
\numberwithin{equation}{section}

\input macros   %abreviations diverses
\begin{document}

\title{Poincar{\'e} lemma and global homotopy formulas with sharp anisotropic H{\"o}lder estimates in q-concave
CR manifolds}

\author{Christine LAURENT-THI\'{E}BAUT}

\date{Pr\'{e}publication de l'Institut Fourier n$^\circ$~xxx (2000)\\
\vspace{5pt} {\tt
http://www-fourier.ujf-grenoble.fr/prepublications.html }}
\date{}
\maketitle

\ufootnote{\hskip-0.6cm  {\it A.M.S. Classification}: 32V20.
}

\bibliographystyle{amsplain}
%\setcounter{page}{-1}
%\input resume

%\tableofcontents
\input an0
\input prelim

\input an1
\input an2

\bibliography{biblio}
\vespa\vespa
\begin{flushright}
  \begin{minipage}[t]{8cm}

Universit\'e de Grenoble

   Institut Fourier

UMR 5582  CNRS/UJF

   BP 74

38402 St Martin d'H\`eres Cedex

France

Christine.Laurent@ujf-grenoble.fr

  \end{minipage}
\end{flushright}

\enddocument

\end

%% file: macros.tex
\providecommand\ufootnote[1]{{\let\thefootnote\relax\footnote[0]{#1}}}

\newcommand{\vespa}{\vspace{1em}}

\newcommand{\ac}{\mathcal A}
\newcommand{\bc}{\mathcal B}

\newcommand{\ic}{\mathcal I}
\newcommand{\cc}{\mathcal C}

\newcommand{\ec}{\mathcal E}

\newcommand{\rb}{\mathbb R}
\newcommand{\nb}{\mathbb N}
\newcommand{\mb}{\mathbb M}
\newcommand{\cb}{\mathbb C}

\newcommand{\ci}{{\mathcal C}^\infty}

\newcommand{\Cal}{\mathcal}

\newcommand{\ol}{\overline}
\newcommand{\ld}{,\ldots,}
\newcommand{\pa}{\partial}
\newcommand{\opa}{\ol\partial}
\newcommand{\wt}{\widetilde}

 \DeclareMathOperator{\im}{Im}
\DeclareMathOperator{\ke}{Ker}

%% file: an0.tex
In this paper we prove sharp anisotropic H{\"o}lder estimates
for the local solutions of the tangential Cauchy-Riemann equation in
$q$-concave CR manifolds and we derive the same kind of estimates for
global solutions when the manifold is compact.

It is wellknown since the works by Folland and Stein \cite{FoSt} that the sharp
H{\"o}lder estimates for the solutions of the $\opa_b$-equation have to be
ansotropic, the complex tangential directions to $M$ playing a special role. 
We consider two types of anisotropic H{\"o}lder spaces the spaces
$\ac^{p+\alpha}$ and the spaces $\Gamma^{p+\alpha}$, which are defined
respectively in Section \ref{s2} and in Section \ref{s3}.

Our first result is a Poincar{\'e} Lemma for the $\opa_b$ operator~:
\begin{thm}\label{poincare}
Let $M$ be a $q$-concave generic CR submanifold of class $\cc^3$ and
of real codimension $k$ of an $n$-dimensional complex manifold $X$ and
$z_0$ a point in $M$. For any open neighborhood $U$ of $z_0$ in $M$,
there exists an open  neighborhood $V\subset U$ of $z_0$ and, for each
$r$ such that $1\leq r\leq q-1$ or $n-k-q+1\leq r\leq n-k$, an operator
$$T_r~:\cc_{n,r}(U)\to\cc_{n,r}(V)$$
with the following properties:

(i)$f_{|_V}=\opa_b T_rf+T_{r+1}\opa_b f$, for $1\leq r\leq q-2$ or
$n-k-q+1\leq r\leq n-k$;

(ii)$f_{|_V}=\opa_b T_rf$, if $r=q-1$ and $\opa_b f=0$;

(iii) if $f\in\ac^{p+\alpha}_{n,r}(U)$, $0<\alpha<1$, then
$T_rf\in\ac^{p+1+\alpha}_{n,r}(V)$ if $M$ is of class
$\cc^{[\frac{p+1}{2}]+3}$ and $1\leq r\leq q-1$ or if $M$ is of class
$\cc^{[\frac{p}{2}]+3}$ and $n-k-q+1\leq r\leq n-k$;

(iv)  if $f\in\Gamma^{p+\alpha}_{n,r}(U)$, $0<\alpha<1$, then
$T_rf\in\Gamma^{p+1+\alpha}_{n,r}(V)$ if $M$ is of class
$\cc^{p+4}$ and $1\leq r\leq q-1$ or if $M$ is of class
$\cc^{p+2}$ and $n-k-q+1\leq r\leq n-k$.
\end{thm}

The operators $T_r$ are the operators defined in \cite{BaLa}. Theorem
\ref{poincare} is already proved in \cite{BaLa} (cf. Theorem 5.9) with $\cc^l$
estimates, it is derived from a local homotopy formula
(cf. Proposition \ref{estmeeslocales} of the present paper).
 The novelty here are the anisotropic H{\"o}lder regularity properties of
the operators $T_r$, which are proved in Sections \ref{s2} and
\ref{s3}.

The second result is a global homotopy formula with sharp anisotropic
H{\"o}lder estimates for compact CR manifolds.
\begin{thm}\label{global}
Let $E$ be an holomorphic vector bundle over $X$ and  $M$ be a
compact, $\ci$-smooth, $q$-concave generic CR submanifold of real codimension
$k$ of an $n$-dimensional complex
manifold $X$. Assume the $\opa_b$-cohomology group $H^{n,r}(M)$
vanishes for some $r$, $0\le r\le q-1$ or $n-k-q+1\le r\le n-k$, then there 
exist  continuous linear
operators
\begin{equation*}
A_s~:~\cc_{n,s}(M,E)\to \cc_{n,s-1}(M,E),~ s=r,r+1
\end{equation*}
 such that:

(i) For all 
$f\in\cc_{n,r}(M,E)$ with $\opa_bf\in\cc_{n,r+1}(M,E)$,
\begin{equation}\label{homotopieglobale}
f=\begin{cases}A_{1}\opa_b f&\text{if }r=0\,,\\\opa_b
A_rf+A_{r+1}\opa_b f\qquad&\text{if }r\geq 1
\end{cases}
\end{equation}

(ii) For all $p\in\nb$ and $0<\alpha<1$ 
\begin{equation*}
A_s(\ac^{p+\alpha}_{n,s}(M,E))\subset \ac^{p+1+\alpha}_{n,s-1}(M,E)
\end{equation*}
and $A_s$ is continuous as an operator between $\ac^{p+\alpha}_{n,s}(M,E)$ and
$\ac^{p+1+\alpha}_{n,s-1}(M,E)$;

(iii) For all $p\in\nb$ and $0<\alpha<1$
\begin{equation*}
A_s(\Gamma^{p+\alpha}_{n,s}(M,E))\subset \Gamma^{p+1+\alpha}_{n,s-1}(M,E)
\end{equation*}
and $A_s$ is continuous as an operator between $\Gamma^{p+\alpha}_{n,s}(M,E)$ and
$\Gamma^{p+1+\alpha}_{n,s-1}(M,E)$.
\end{thm}

The operators $A_r$ are derived from the local operators $T_r$ by the
globalization process from \cite{LaLeglobal} and \cite{BroLaLe}. Then
Theorem \ref{global} follows immediately from Theorem
\ref{formuleglobale} and the estimates in Sections \ref{s2} and
\ref{s3}. 

The Poincar{\'e} Lemma  for the $\opa_b$ operator with sharp anisotropic
H{\"o}lder estimates is new in the case of CR manifolds of arbitrary
codimension, 
it was proved first for the Heisenberg group and more generally for
strictly pseudoconvex hypersurfaces by Folland and Stein \cite{FoSt}.
A related result in CR manifolds of arbitrary
codimension is a local almost homotopy formula with sharp anisotropic
H{\"o}lder estimates proved by Polyakov in \cite{Polyanis}, but there is a
mistake in the construction of the barriers which is corrected in
\cite{Poly1} and we hope that the same anisotropic estimates hold for
his new kernels.

The global result is not totally new,
it is proved in \cite{ShWa} even for abstract CR manifolds for the
anisotropic Sobolev spaces (Sobolev version of the $\ac^{p+\alpha}$
spaces) and for the Folland-Stein spaces 
$\Gamma^{p+\alpha}$ if $r\neq q-1$. The proof is based on $L^2$
theory for the $\Box_b$ operator and the Hodge decomposition theorem
which gives a global homotopy formula in degree $r$ but only if the CR manifold
$M$ is $(r+2)$-concave.

Our result contains the important case where $M$ is only $2$-concave
and $r=1$, which may be useful in the study of the embeddability of
compact CR manifolds (in \cite{ShWa}, results in degree $1$ need the
manifold to be $3$-concave). 

%%% Local Variables: 
%%% mode: latex
%%% TeX-master: "anisotropic"
%%% End: 

%% file: prelim.tex
\section{Preliminaries and definitions}\label{s1}

Let $(\mb,H_{0,1}\mb)$ be a generically embeddable abstract compact CR manifold of class
$\ci$ and $\ec~:~\mb\to M\subset X$ be a
$\ci$-smooth CR generic embedding of $\mb$ in a complex manifold $X$, then
$M$ is a compact CR
submanifold of $X$ of class $\ci$ with the CR structure
$H_{0,1}M=d\ec(H_{0,1}\mb)=T_\cb M\cap T_{0,1} X$ and
the tangential Cauchy-Riemann operator $\opa_b$ induced by the
Cauchy-Riemann operator $\opa$ from the complex manifold $X$.

An generic CR manifold $\mb$ is said to be $q$-concave if its
Levi form at each point admits at least $q$ negative eigenvalues in
all directions.

Assume $\mb$ is $q$-concave, then $M$ is also $q$-concave and we may
apply the results in \cite{BaLa} and \cite{LaLeglobal},
\cite{BroLaLe} on local estimates and global homotopy formulas for
the tangential Cauchy-Riemann operator.

First let us recall the definition of the usual H\"older spaces of
forms. If $D$ is relatively compact domain in $X$, then

- $\cc^\alpha(\ol D\cap M)$, $0\leq\alpha<1$, is the set of
continuous functions on $\ol D\cap M$ which are
H\"older continuous with exponent $\alpha$ on $\ol D\cap M$ , if
$\alpha>0$. We set 
\begin{equation}\label{holder}
\|f\|_\alpha=\sup_{z\in D\cap M}|f(z)|+\sup_{z,\zeta\in
  D\cap M\atop z\neq\zeta}\frac{|f(z)-f(\zeta)|}{|z-\zeta|^\alpha}
\end{equation}

- $\cc^{l+\alpha}(\ol D\cap M)$, $l\in\nb$,
$0\leq\alpha<1$, is the space of functions of
class $\cc^l$ on $\ol D\cap M$, whose derivatives of order $l$ are in
$\cc^\alpha(\ol D\cap M)$.

The H\"older space $\cc^{l+\alpha}_*(\ol D\cap M)$, $l\in\nb$,
$0\leq\alpha<1$, is
then the space of continuous forms on $\ol D\cap M$, whose
coefficients are in $\cc^{l+\alpha}(\ol D\cap M)$.

In \cite{BaLa} the following result is proved

\begin{prop}\label{estmeeslocales}
Let $M$ be a $q$-concave generic CR
submanifold of $X$ of class $\ci$.
For each point in $M$, there exist
a neighborhood $U$ and linear operators
\begin{equation*}
T_r~:~\cc^0_{n,r}(M)\to \cc^0_{n,r-1}(U)\,,\quad 1\leq r\leq q~ {\rm
and}~ n-k-q+1\le r\le n-k\,,
\end{equation*}
with the following two properties~:

(i) For all $l\in\nb$ and $1\leq r\leq q$ or $n-k-q+1\le r\le n-k$,
\begin{equation*}
T_r(\cc^l_{n,r}(M))\subset \cc^{l+1/2}_{n,r-1}(\ol U)
\end{equation*}
and $T_r$ is continuous as an operator between $\cc^l_{n,r}(M)$ and
$\cc^{l+1/2}_{n,r-1}(\ol U)$.

(ii) If $f\in\cc^1_{n,r}(M)$, $0\leq r\leq q-1$ or $n-k-q+1\le r\le
n-k$, has compact support in $U$, then on $U$,
\begin{equation}\label{homotopielocale}
f=\begin{cases}T_{1}\opa_b f&\text{if }r=0\,,\\\opa_b
T_rf+T_{r+1}\opa_b f\qquad&\text{if }1\leq r\leq q-1 ~ {\rm or}~
n-k-q+1\le r\le n-k \,.
\end{cases}
\end{equation}
\end{prop}
{\parindent=0pt and in \cite{LaLeglobal} and \cite{BroLaLe} we have
derive from the
  previous proposition
a global homotopy formula by mean of a functional analytic
construction:}

\begin{thm}\label{formuleglobale}
Let $E$ be an holomorphic vector bundle over $X$ and  $M$ be a
compact $q$-concave generic CR submanifold of $X$ of class $\ci$. Then there
exist finite dimensional subspaces $\Cal H_r$ of $\Cal
Z^\infty_{n,r}(M,E)$, $0\le r\le q-1$ and $n-k-q+1\le r\le n-k$,
where $\Cal H_0=\Cal Z^\infty_{n,0}(M,E)$, continuous linear
operators
\begin{equation*}
A_r~:~\cc^0_{n,r}(M,E)\to \cc^0_{n,r-1}(M,E),\quad 1\leq r\leq q~
{\rm and}~ n-k-q+1\le r\le n-k
\end{equation*}
and continuous linear projections $$ P_r:\Cal
C^0_{n,r}(M,E)\rightarrow \Cal C^0_{n,r}(M,E) \,,\qquad 0\le r\le
q-1~ {\rm and}~ n-k-q+1\le r\le n-k\,,
$$with
\begin{equation}\label{19.4.07'}
\im P_r=\Cal H_r\,,\qquad 0\le r\le q-1~ {\rm and}~ n-k-q+1\le r\le
n-k\,,
\end{equation}and
\begin{equation}\label{20.4.07}
\cc^0_{n,r}(M,E)\cap \opa_b\cc^0_{n,r-1}(M,E)\subseteq\ke P_r\,,\qquad 1\le
r\le q-1~ {\rm and}~ n-k-q+1\le r\le n-k,
\end{equation} such that:

(i) For all $l\in\nb$ and $1\leq r\leq q$ or $n-k-q+1\le r\le n-k$,
\begin{equation*}
A_r(\cc^l_{n,r}(M,E))\subset \cc^{l+1/2}_{n,r-1}(M,E)
\end{equation*}
and $A_r$ is continuous as an operator between $\cc^l_{n,r}(M,E)$ and
$\cc^{l+1/2}_{n,r-1}(M,E)$.

(ii) For all $0\leq r\leq q-1$ or $n-k-q+1\le r\le n-k$ and
$f\in\cc^0_{n,r}(M,E)$ with $\opa_bf\in\cc^0_{n,r+1}(M,E)$,
\begin{equation}\label{homotopieglobale}
f-P_rf=\begin{cases}A_{1}\opa_b f&\text{if }r=0\,,\\\opa_b
A_rf+A_{r+1}\opa_b f\qquad&\text{if }1\leq r\leq q-1 ~ {\rm or}~
n-k-q+1\le r\le n-k \,.
\end{cases}
\end{equation}
\end{thm}

We have now to recall the main steps of the construction of the local
kernels defining the operators $T_r$.

Let $M$ be a generic $CR$ manifold of class $\cc^3$ in $\cb^n$,
$z_0$ a point in $M$ and $U_0$ an open neighborhood of $z_0$ in
$\cb^n$ and $\widehat\rho_1\ld \widehat\rho_k$ some functions of
class $\cc^3$ from $U_0$ into $\rb$ such that
$$M\cap U_0 = \{z\in U_0 ~|~\widehat\rho_1 (z)=\dots=\widehat\rho_k
(z)=0\}$$ and satisfying $\opa\widehat\rho_1
(z)\wedge\dots\wedge\opa\widehat\rho_k (z)\neq 0$ for $z\in M\cap
U_0$.

Let $C>0$ be a fixed constant, we set, for $j=1\ld k$,
\begin{equation}\label{equ5}
\begin{aligned}
   &\rho_j = \widehat\rho_j + C\sum_{\nu=1}^k\widehat\rho_{\nu}^2\\
  &\rho_{-j} = -\widehat\rho_j +
  C\sum_{\nu=1}^k\widehat\rho_{\nu}^2.
\end{aligned}
\end{equation}

We define $\ic$ as the set of all subsets $I\subset\{\pm 1\ld\pm
k\}$ such that $|i|\neq |j|$ for all $i,j\in I$ with $i\neq j$.
For $I\in \ic$, $|I|$ denotes the number of elements in $I$, then
$\ic(l)$, $1\leq l\leq k$, is the set of all $I\in \ic$ with
$|I|=l$ and $\ic'(l)$, $1\leq l\leq k$, is the set of all $I\in
\ic$ of the form $I=(i_1\ld i_l)$ with $|i_{\nu}|=\nu$ for
$\nu=1\ld l$.

If $I\in \ic$ and $\nu\in\{1\ld |I|\}$, then $i_\nu$ is the
element with rank $\nu$ in $I$ after ordering $I$ by modulus.
We set $I(\widehat\nu)=I\setminus\{i_\nu\}$.

If $I\in \ic$, then

\hskip2cm{${\rm sgn}I=1$ if the number of negative elements in
$I$ is even}

\hskip2cm{${\rm sgn}I=-1$ if the number of negative elements in
$I$ is odd.}

Let $(e_1\ld e_k)$ be the canonical basis of $\rb^k$, set
$e_{-j}=-e_j$ for every $1\leq j\leq k$. Let $I=(i_1\ld i_l)$ be
in $\ic(l)$, $1\leq l\leq k$, set
$$\wt\Delta_I=\{x=\sum_{j=1}^l \lambda_j e_{i_j}|\lambda_i\geq 0,
1\leq i\leq l, \sum_{i=1}^l \lambda_i=1\}.$$
We identify the abstract simplex $\Delta_I$ with the geometric simplex
$\wt\Delta_I$ by setting $x(\lambda)=\sum_{j=1}^l \lambda_j e_{i_j}$
for all $\lambda\in\Delta_I$

For all $I\in\ic'(k)$, we denote by $I*$ the multi-index $(i_1\ld
i_{k},*)$, where $I=(i_1\ld i_{k})$, and by $\ic'(k,*)$ the set of
all multi-indexes $I*$, with $I\in\ic'(k)$. We set
$\rho_*=\frac{1}{k}(\rho_1+\dots +\rho_k)$  and $\rho_\lambda
=\lambda_1\rho_{i_1}+\dots +\lambda_k\rho_{i_k} +\lambda_*\rho_*$
for $\lambda=(\lambda_1\ld \lambda_k,\lambda_*)\in \Delta_{I*}$.

We denote by $D$ a relatively compact open subset of $U_0$  and for
$I\in\ic$, $I=(i_1\ld i_{|I|})$, we define
\begin{align*}
  & D_I=\{\rho_{i_1}<0\}\cap\dots\cap\{\rho_{i_{|I|}}<0\}\cap D\\
  & D_I^*=\{\rho_{i_1}>0\}\cap\dots\cap\{\rho_{i_{|I|}}>0\}\cap D\\
  & S_I=\{\rho_{i_1}=0,\dots ,\rho_{i_{|I|}}=0\}\cap D\\
  & \Gamma_I=\{\rho_{i_1}=\dots =\rho_{i_{|I|}}\}\cap D_I\\
  & \Gamma_I^*=\{\rho_{i_1}=\dots =\rho_{i_{|I|}}\}\cap D_I^*.
\end{align*}
These manifolds are oriented as follows~: $D_I$ and $D^*_I$ as
$\cb^n$ for all $I\in\ic$, $S_{\{j\}}$ as the boundary of
$D_{\{j\}}$ for $j=\pm 1\ld \pm k$, $S_I$ as the boundary of
$S_{I(\widehat {|I|})}\cap \ol D_{\{i_{|I|}\}}$ for all $I\in\ic$,
$|I|\geq 2$, $\Gamma_I$ such that $S_I=\pa\Gamma_I$ and $M\cap D$ as
$S_I$ with $I=\{1\ld k\}$.

If $M$ is $q$-concave, it follows from Lemma 3.1.1 in~\cite{AiHe}
that we can choose the constant $C$ in (\ref{equ5}) such that the
functions $\rho_j$, $-k\leq j\leq k$, $j\neq 0$, have the following
property: for each $I\in\ic'(k)$ and every $\lambda\in\Delta_I$, the
Levi form of the defining function $\rho_\lambda$ of $M$ in the
direction $x(\lambda)$ has at least $q+k$ positive eigenvalues on
$U'\subset\subset U_0$. Then using the method developed in section 3
of \cite{LaLefourier} and the results in \cite{LaLeregular}, we can
construct, for each $\lambda\in\Delta_I$ a Leray section in the
direction $x(\lambda)$, which has some holomorphy properties and
depends smoothly on $\lambda$ and we get on $U'\times U'\setminus
\Delta(U')$, $U'\subset\subset U_0$, some Cauchy-Fantappi\'e kernels
$$C_{I*}(z,\zeta)=\int_{\lambda\in\Delta_{I*}}K_{I*}(z,\zeta,\lambda)$$
for each $I\in\ic'(k)$ (cf.\cite{BaLa}) such that, if we set
$R_{M}=\sum_{I\in\ic'(k)} {\rm sgn}(I) C_{I*}$, the associated
integral operators $T_r$ satisfy the homotopy formula (ii) of
Proposition \ref{estmeeslocales} with $U=D\cap M$.

We can describe the singularity of the kernels $C_{I*}$  in the
following way.

A form of type $O_s$ (or of type $O_s(z,\zeta,\lambda)$) on $\ol
D_I\times \ol D^*_I\times \Delta_{I*}$ is, by definition, a
continuous differential form $f(z,\zeta,\lambda)$ defined for all
$(z,\zeta,\lambda)\in\ol D_I^*\times \ol D_I\times
\Delta_{I*}$ with $z\neq \zeta$ such that the following
conditions are fulfilled~:
\begin{enumerate}
  \item All derivatives of the coefficients of $f$ which are of order $0$ in
 $z$, and of order $\leq 1$ in $\zeta$ and of arbitrary order in
 $\lambda$ are continuous
for all $(z,\zeta,\lambda)\in\ol D_I^*\times \ol D_I\times
\Delta_{I*}$ with $z\neq \zeta$.
  \item  Let $\nabla^\kappa_\zeta$, $\kappa=0,1$, be a
  differential operator with constant
  coefficients, which is of order $0$ in $z$, of order $\kappa$ in
   $\zeta$ and of arbitrary order in $\lambda$. Then there is a constant
  $C>0$ such that, for each coefficient $\varphi(z,\zeta,\lambda)$ of the form
  $f(z,\zeta,\lambda)$,
  \begin{equation*}
|\nabla^\kappa_\zeta\varphi(z,\zeta,\lambda)|\leq
C|\zeta-z|^{s-\kappa}
  \end{equation*}
for all $(z,\zeta,\lambda)\in\ol D_I\times \ol D_I^*\times
\Delta_{I*}$ with $z\neq \zeta$.
\end{enumerate}

Following the calculations in \cite{BaLa}, we get
\begin{equation}\label{singRC}
[K_{I*}(z,\zeta,\lambda )]_{\rm{deg}\lambda=|I|}\wedge
dz_1\wedge\dots\wedge dz_n=\sum_{0\leq m\leq k\atop i_1\ld i_m\in
I}\frac{O_{|I|+1-m}}{\Phi^n}\wedge \pa\rho_{i_1}(z)\wedge\dots\wedge
\pa\rho_{i_m}(z).
\end{equation}
and
\begin{equation}\label{singdbarRC}
[\opa_z K_{I*}(z,\zeta,\lambda )]_{\rm{deg}\lambda=|I|}\wedge
dz_1\wedge\dots\wedge dz_n=\sum_{0\leq m\leq k\atop i_1\ld i_m\in
I}\frac{O_{|I|-m}}{\Phi^n}\wedge \pa\rho_{i_1}(z)\wedge\dots\wedge
\pa\rho_{i_m}(z).
\end{equation}

The support function associated to the Leray section used to construct
the kernels $C_{I*}$ satisfies for $\zeta,z$ in a neighborhood of
$U'$ and $\lambda\in\Delta_{I*}$
\begin{equation}\label{phi}
\Phi(z,\zeta,\lambda)=2\sum_{j=1}^k\frac{\pa\rho_\lambda}{\pa\zeta_j}(\zeta)(\zeta_j-z_j)+O_2
\end{equation}
and
\begin{equation}\label{Rephi}
{\rm Re}~\Phi(z,\zeta,\lambda)\geq\rho_\lambda(\zeta)-\rho_\lambda(z)+\gamma|\zeta-z|^2.
\end{equation}
Consequently, if $X_\cb$ denotes a complex tangent vector field to
$M$, then
\begin{equation}\label{phitan}
X^\zeta_\cb\Phi(z,\zeta,\lambda)=O_1 \quad {\rm and}\quad
X^z_\cb\Phi(z,\zeta,\lambda)=O_1.
\end{equation}
Moreover the function $\Phi(z,\zeta,\lambda)$ is of class $\ci$ in its
first variable $z$
and of class $\cc^{l-1}$ in its second variable $\zeta$ if the
manifold $M$ is of class $\cc^l$.

In the next sections we shall prove that in fact the operators $T_r$
and $A_r$, $1\leq r\leq q$, satisfy sharp anisotropic estimates.

%% file: an1.tex
\section{First anisotropic estimates}\label{s2}

If $M$ is a generic CR submanifold of a complex manifold $X$, the
 tangent bundle to $M$ admits a maximal complex subbundle 
called the complex tangent bundle to $M$ and denoted by $HM$. It is
related to the CR structure of $M$ by $HM=TM\cap (H_{0,1}M+\ol H_{0,1}M)$.

Let us define now some anisotropic H\"older spaces of forms which
specify the complex tangent directions:

- $\ac^\alpha(\ol D\cap M)$, $0<\alpha<1$, is the
set of continuous functions on $\ol D\cap M$ which are in
$\cc^{\alpha/2}(\ol D\cap M)$.

- $\ac^{1+\alpha}(\ol D\cap M)$, $0<\alpha<1$, is the
set of functions $f$ such that $f\in\cc^{(1+\alpha)/2}(\ol D\cap M)$ and
$X_\cb f\in\cc^{\alpha/2}(\ol D\cap M)$, for all
complex tangent vector fields $X_\cb$ to $M$.  Set
\begin{equation}\label{holderanis}
  \|f\|_{A\alpha}=\|f\|_{(1+\alpha)/2}+\sup_{\|X_\cb\|\leq 1}\|X_\cb f\|_{\alpha/2}
\end{equation}

- $\ac^{p+\alpha}(\ol D\cap M)$, $p\geq 2$, $0<\alpha<1$, is the
set of functions $f$ of class $\cc^{[p/2]}$ such that
$Xf\in\ac^{p-2+\alpha}(\ol D\cap M)$,
for all tangent vector fields $X$ to $M$ and
$X_\cb f\in\ac^{p-1+\alpha}(\ol D\cap M)$, for all
complex tangent vector fields $X_\cb$ to $M$.

The anisotropic H\"older space of forms $\ac^{p+\alpha}_*(\ol D\cap
M)$, $p\geq 0$, $0<\alpha<1$, is then the space of continuous forms on
$\ol D\cap M$, whose
coefficients are in  $\ac^{p+\alpha}(\ol D\cap M)$.

This section is devoted to the study of the continuity properties
with respect to 
these spaces of
the operators $T_r$
and $A_r$, $1\leq r\leq q$ and $n-k-q+1\leq r\leq n-k$, defined in
the previous section when $M$ 
is $q$-concave.

We first prove some estimates for the operators defined by the local
kernels from the previous section. For this we may assume that
$M$ is embedded in $\cb^n$.
To clarify the exposition, we introduce a new kind of operators.

\begin{defin}
Let $m\in\nb$, $\epsilon\in\{-1,0,1\}$, $\delta>0$ and $I\in\ic'(k)$. An
\emph{operator of type} $(m,\epsilon,\delta)$ is, by definition, a map
$$E~:~\cc^0_{n,*}(\Gamma_I)\to\cc^\infty_{n,*}(\Gamma_I)$$
such that there exist
\begin{itemize}
\item an integer $\kappa\geq 0$
\item a differential form $\widehat E(z,\zeta,\lambda)$ of type
  $O_{|I|+1-2n+2\kappa+m+\epsilon}$ on $\Gamma_I\times
  D^*_I\times\Delta_{I*}$ such that for all $f\in\cc^0_{n,*}(\Gamma_I)$
$$Ef(\zeta)=\int_{(z,\lambda)\in\Gamma_I\times\Delta_{I*}}
\wt f(z)\wedge\frac{\widehat
  E(z,\zeta,\lambda)\wedge\Theta(z)}{(\Phi+\delta)^{\kappa+m}(z,\zeta,\lambda)},$$
where$\wt f\in \cc^0_{0,*}(\Gamma_I)$ is the form with
$$f(z)=\wt f(z)\wedge dz_1\wedge\dots\wedge dz_n,$$
and $\Theta=1$, if $m=0$, and either
$\Theta=\pa\rho_{i_1}\wedge\dots\wedge\pa\rho_{i_m}$ or
$\Theta=\opa\rho_{i_1}\wedge\dots\wedge\pa\rho_{i_m}$, if $m\geq 1$
with $i_1,\dots,i_m\in I$.
\end{itemize}
\end{defin}

We consider the operators $\wt C_{I*}$, $I\in\ic'(k)$, defined by
\begin{equation}\label{operateurC}
\wt C_{I*} f(\zeta)=\int_{z\in D\cap M}f(z)\wedge
C_{I*}(z,\zeta)=\int_{(z,\lambda)\in D\cap
  M\times\Delta_{I*}}f(z)\wedge  K_{I*}(z,\zeta,\lambda)
\end{equation}
for $f\in\cc^1_{n,r}(D)$, $0\leq r\leq n-k$.
Using Stokes formula we get
\begin{equation}
\wt C_{I*} f(\zeta)=\int_{(z,\lambda)\in\Gamma_I\times\Delta_{I*}}\opa
f(z)\wedge
K_{I*}(z,\zeta,\lambda)+(-1)^{n+r}\int_{(z,\lambda)\in\Gamma_I\times\Delta_{I*}}f(z)\wedge
\opa_z K_{I*}(z,\zeta,\lambda).
\end{equation}

By (\ref{singRC}) and ((\ref{singdbarRC}) the kernel $\wt C_{I*}$ is
the sum of an operator of type 
$(m,1,0)$ and of an operator of type $(m,0,0)$, $0\leq m\leq k$, and
if $X_\cb$ denotes a complex tangent vector field to $M$, then by
(\ref{phitan}) the kernels 
$X^\zeta_\cb\wt C_{I*}$ and $X^z_\cb\wt C_{I*}$ are the sum of an
operator of type $(m,0,0)$ and of an operator of type $(m,-1,0)$,
$0\leq m\leq k$.

To prove the anisotropic estimates we need the following spaces
and norms on differential forms~:

- $\bc^\beta_*(\Gamma_I)$, $\beta\geq 0$, is the space of forms
$f\in\cc^0_*(\ol\Gamma_I\setminus M)$ such that, for some constant $C>0$,
$$\|f(z)\|\leq C[{\rm dist}(z,M)]^{-\beta}, \quad z\in\Gamma_I,$$ where
${\rm dist}(z,M)$ is the Euclidean distance between $z$ and $M$. Set
\begin{equation}\label{croissance}
\|f\|_{\beta}=\sup_{z\in \Gamma_I}\big(\|f(z)\|[{\rm
dist}(z,M)]^{\beta}\big)
\end{equation}
for $\beta\geq 0$ and $f\in\bc^\beta_*(\Gamma_I)$.

It follows from the characterization of H\"older continuous functions
by mean of the Poisson integral that, if $f\in\cc^\alpha_*(\ol D\cap M)$,
$0<\alpha<1$, there exists for each $I\in\ic'(k)$ a continuous form
$f_I$ on $\ol\Gamma_I$ such that ${f_I}_{|_{M}}=f$ and
$X_If_I\in\bc^{1-\alpha}_*(\Gamma_I)$, for all vector
fields $X_I$ tangent to $\Gamma_I$.

\begin{lem}\label{estimbeta}
If $f\in \bc^\beta_*(\Gamma_I)$, $0\leq \beta<1$, has compact
support in $D$ and if $E_I$ is an operator of type
$(m,\epsilon,\delta)$, then there exists a constant $C>0$,
independent of $\delta$, such that for $\zeta\in D^*_I$
\begin{align*}
|E_I f(\zeta)|&\leq C \|f\|_{\beta}[{\rm
  dist}(\zeta,M)+\delta]^{1/2-\beta},\quad {\rm if}~ \epsilon=0\\
|E_I f(\zeta)|&\leq C \|f\|_{\beta}[{\rm
  dist}(\zeta,M)+\delta]^{-\beta},\quad {\rm if}~ \epsilon=-1\\
|E_I f(\zeta)|&\leq C \|f\|_{\beta}[{\rm
  dist}(\zeta,M)+\delta]^{1-\beta},\quad {\rm if}~ \epsilon=1
\end{align*}
\end{lem}
\begin{proof}
It follows from \cite{LaLefourier} and \cite{BaLa} that, after
integration in $\lambda$, an operator of type $(m,\epsilon,\delta)$ is
controled, since $n\geq 3$, by
\begin{equation*}
\int_{\Gamma_I}\frac{|f_I|~|\sigma\wedge d\rho_I{\wedge_{\nu=1}^s}dt_\nu|}
{(|\rho_I|+d+\delta+|\zeta-z|^2)\Pi_{\nu=1}^{s}(|t_\nu|+d+\delta+|\zeta-z|^2)^{1+1/s} 
  |\zeta-z|^{2n-k-s-3-\epsilon}},
\end{equation*}
with $1\leq s\leq k$,
 where $\rho_I$ is the function defined by
$\rho_I(z)=\rho_{i_1}(z)=\dots=\rho_{i_k}(z)$ for $z\in\Gamma_I$,
$d=d(\zeta)={\rm dist}(\zeta,M)$,
$\sigma$ a monomial in $dz_1,\dots,dz_n,d\ol z_1,\dots,d\ol z_n$, $t_\nu= {\rm
  Im}~\Phi(z,\zeta,\lambda^\nu)$ and $dt_\nu= d_z{\rm
  Im}~\Phi(z,\zeta,\lambda^\nu)$ with $\lambda^1,\dots,\lambda^{k+1}$
some points in $\Delta_{I*}$ which are 
linearly independent in $\rb^{k+1}$,

We denote all the constants by $C$. Then since $d(z)\leq
C|\rho_I(z)|$ for $z\in\Gamma_I$, we get
\begin{align*}
\|E_I &f(\zeta)\|\leq\\
& C
\|f\|_{\beta}\int_{\Gamma_I}\frac{|\rho_I|^{-\beta} ~|\sigma\wedge
  d\rho_I{\wedge_{\nu=1}^s}dt_\nu|}
{(|\rho_I|+d+\delta+|\zeta-z|^2)\Pi_{\nu=1}^{s}(|t_\nu|+d+\delta+|\zeta-z|^2)^{1+1/s}
  |\zeta-z|^{2n-k-s-3-\epsilon}}
\end{align*}
for all $f\in \bc^\beta_*(\Gamma_I)$ and $\zeta\in D^*_I$. We can
use $\rho_I$ as a coordinate in $\Gamma_I$ and since $M$ is generic
the $t_\nu$'s can also be used as coordinates, so we obtain that
\begin{align*}
\|E_I &f(\zeta)\|\leq \\
&C
\|f\|_{\beta}\int_{y\in\rb^{2n-k+1}\atop |y|<c}\frac{|y_1|^{-\beta}
    ~dy_1\wedge dy_2\wedge\dots
    \wedge dy_{2n-k+1}}
{(|y_1|+d+\delta+|\zeta-z|^2)\Pi_{\nu=1}^{s}(|y_\nu|+d+\delta+|y|^2)^{1+1/s}
  |y|^{2n-k-s-3-\epsilon}}.
\end{align*}
Since $\beta<1$, we can integrate with respect to $y_1$, which gives
\begin{equation*}
\|E_I f(\zeta)\|\leq C
\|f\|_{\beta}\int_{y\in\rb^{2n-k}\atop |y|<c}\frac{dy_2\wedge\dots
    \wedge dy_{2n-k+1}}
{(d+\delta+|y|^2)^{\beta}\Pi_{\nu=1}^{s}(|y_\nu|+d+\delta+|y|^2)^{1+1/s}
  |y|^{2n-k-s-3-\epsilon}}.
\end{equation*}
Then integrating with respect to $y_2,\dots,y_s$ and using spherical coordinates,
we have
\begin{equation*}
\|E_I f(\zeta)\|\leq C
\|f\|_{\beta}\int_0^C\frac{dr}
{(d+\delta+r^2)^{\beta+1} r^{-2-\epsilon}},
\end{equation*}
which proves the lemma following the values of $\epsilon$.
\end{proof}

In the same way, following again \cite{LaLefourier}, we can prove
estimates for the gradient of $E_If$, when  $E_I$ is an
operator of type $(m,\epsilon,\delta)$ and $f\in \bc^\beta_*(\Gamma_I)$,
$0\leq \beta<1$.

\begin{lem}\label{estimbeta-gradient}
If $f\in \bc^\beta_*(\Gamma_I)$, $0\leq \beta<1$, has compact support
in $D$,  $E_I$ is  an
operator of type $(m,\epsilon,\delta)$ and if $\nabla_\zeta$ is one of the
operators $\frac{\pa}{\pa\zeta_1},\dots,\frac{\pa}{\pa\zeta_n}$ or
$\frac{\pa}{\pa\ol\zeta_1},\dots,\frac{\pa}{\pa\ol\zeta_n}$, then
there exists a constant $C>0$ independent of $\delta$ such
that for $\zeta\in D^*_I$
\begin{align*}
|\nabla_\zeta E_I f(\zeta)|&\leq C \|f\|_{\beta}[{\rm
  dist}(\zeta,M)+\delta]^{-1/2-\beta},\quad {\rm if}~ \epsilon=0\\
|\nabla_\zeta E_I f(\zeta)|&\leq C \|f\|_{\beta}[{\rm
  dist}(\zeta,M)+\delta]^{-1-\beta},\quad {\rm if}~ \epsilon=-1\\
|\nabla_\zeta E_I f(\zeta)|&\leq C \|f\|_{\beta}[{\rm
  dist}(\zeta,M)+\delta]^{-\beta},\quad {\rm if}~ \epsilon=1
\end{align*}
\end{lem}

Using the classical Hardy-Littlewood lemma, we deduce from Lemma
\ref{estimbeta}
and Lemma \ref{estimbeta-gradient}, the following estimates for the
operators $\wt C_{I*}$ and $X^\zeta_\cb\wt C_{I*}$

\begin{prop}\label{estimci}
Let $f\in \bc^\beta_*(\Gamma_I)$, $0\leq \beta<1$, be a form with
compact support in $D$, then
\begin{align*}
&\wt C_{I*}f\in\cc^{1/2-\beta}_*(D^*_I),\quad {\rm if}~ 0\leq\beta<1/2\\
&\wt C_{I*}f\in\bc^{\beta-1/2}_*(D^*_I),\quad {\rm if}~ 1/2\leq\beta<1
\end{align*}
and
\begin{equation*}
X^\zeta_\cb\wt C_{I*}f\in \bc^\beta_*(D^*_I).
\end{equation*}
\end{prop}

Let us recall Lemma 5.3 and 5.5 from \cite{BaLa}

\begin{lem}
There exists $Y_1,\dots,Y_k$ tangential vector fields to $M$ such that
for all $\zeta\in U_0$ and all $1\leq i,j\leq k$
$$Y_i^\zeta\Phi_j(\zeta,\zeta)=\delta_{ij},$$
where $\delta_{ij}$ denotes the Kronecker index.
\end{lem}

If $X^z$ is a vector field on $D$ in the variable $z$, we denote by
$X^\zeta$ the same vector field in the variable $\zeta$.

\begin{lem}\label{derivation}
Let us fix $\delta>0$, if $\widehat E_\delta$ is the kernel of an
operator of type $(m,\epsilon,\delta)$, there exists $\eta>o$ such
that for $(z,\zeta)\in \ol D_I\times\ol D_I^*$ with $|z-\zeta|<\eta$
and $z\neq \zeta$ we have
\begin{equation}\label{eq-derivation}
X^z\widehat E_\delta=-X^\zeta\widehat
E_\delta+\frac{(X^z+X^\zeta)\Phi_{I*}}{Y^\zeta_\lambda\Phi_{I*}}Y^\zeta_\lambda
\widehat E_\delta+\widehat G_\delta,
\end{equation}
where $\widehat G_\delta$ is a finite sum of kernels associated to
operators of  type $(m,\epsilon,\delta)$.
\end{lem}

As regularity is a local problem, let us fix $\zeta_0$ in $D\cap M$
and choose a function $\chi$ with compact support  such that
$\chi(z)=1$, if $|z-\zeta_0|\leq \eta/4$, and $\chi(z)=0$, if
$|z-\zeta_0|\geq \eta/2$, where $\eta$ is given by Lemma
\ref{derivation}.

If $f$ is a continuous form on $D\cap M$ and $f_I$ a continuous
extension of $f$ to $\Gamma_I$, we write
\begin{equation*}
\wt C_{I*} f(\zeta)=\wt C_{I*} (\chi f)(\zeta)+\wt C_{I*} ((1-\chi)f)(\zeta)
\end{equation*}
As the singularity of the kernel associated to the operator $\wt
C_{I*}$ is concentrated on the diagonal, the regularity of the second
term, for $\zeta\in M$ with $|\zeta-\zeta_0|\leq \eta/4$, is directly
given by the  regularity of this kernel, which is of class $\ci$ in
the first variable $z$ and of
class $\cc^{l-2}$ in the second variable $\zeta$ when $M$ is of class $\cc^l$.
\begin{thm}\label{estimlocale}
Assume $M$ is of class $\cc^l$.
The operator $\wt R_{M}=\sum_{I\in\ic'(k)} {\rm sgn}(I)
\wt C_{I*}$ is a bounded linear operator from $\ac^{p+\alpha}_*(\ol D\cap
M)$ into $\ac^{p+1+\alpha}_*(\ol D\cap
M)$ for  $0\leq p\leq 2l-7$ and $0<\alpha<1$.
\end{thm}

\begin{proof}
Let us first consider the case $p=0$. Let $f\in\ac^{\alpha}_*(\ol D\cap
M)$, i.e. $f\in\cc^{\alpha/2}_*(\ol D\cap M)$, then there exists, for each
$I\in\ic'(k)$, a continuous form
$f_I$ on $\ol\Gamma_I$ such that ${f_I}_{|_{M}}=f$ and
$X_If_I\in\bc^{1-\alpha/2}_*(\Gamma_I)$, for all vector
fields $X_I$ tangent to $\Gamma_I$. From the previous remark we have
only to study $\wt C_{I*} (\chi f)(\zeta)$.
Let us denote by $\wt C_{I*}^\delta$ the operator $\wt C_{I*}$, in
which we have replaced
$\Phi$ by $\Phi+\delta$. From Lemma
\ref{derivation}, we get that for all vector fields $X_I$ tangent to
$\ol\Gamma_I\cup\ol\Gamma_I^*$, we
have $X_I^\zeta\wt C_{I*}^\delta \chi f=\wt
C_{I*}^\delta X_I^z \chi f_I+E_I \chi f$, where $E_I$ is a sum of operators of
type $(m,\epsilon,\delta)$, $\epsilon\geq 0$.
Let also $X_\cb$ be a complex tangential vector field to $M$.

We deduce from Lemma \ref{estimbeta} that, for $\zeta\in \Gamma_I^*$,
$$|X_I^\zeta \wt C_{I*}^\delta
(\chi f)(\zeta) |\leq C~ \|X_I(\chi f_I)\|_{\beta}[{\rm
  dist}(\zeta,M)+\delta]^{(\alpha-1)/2}$$
and
$$|X_I^\zeta X_\cb\wt C_{I*}^\delta
(\chi f)(\zeta) |\leq C'~ \|X_If_I\|_{\beta}[{\rm
  dist}(\zeta,M)+\delta]^{\alpha/2-1},$$
 where $C$ and $C'$ are constants
independent of $\delta$. Moreover $X_I^\zeta \wt C_{I*}^\delta f$
and $X_I^\zeta X_\cb\wt C_{I*}^\delta f$ converge respectively to
$X_I^\zeta \wt C_{I*} f$ and $X_I^\zeta X_\cb\wt C_{I*} f$ uniformly
on each compact subset of $\Gamma_I^*$, when $\delta$ tends to zero.
This implies that $X_I^\zeta \wt C_{I*}
f\in\bc^{(1-\alpha)/2}_*(\Gamma_I^*)$ and $X_I^\zeta X_\cb\wt C_{I*}
f\in\bc^{1-\alpha/2}_*(\Gamma_I^*)$ and by the classical
Hardy-Littlewood lemma that $\wt C_{I*} f\in\cc^{(\alpha+1)/2}_*(\ol
D\cap M)$ and $X_\cb\wt C_{I*} f\in\cc^{\alpha/2}_*(\ol D\cap M)$,
which means that $\wt C_{I*} f\in\ac^{1+\alpha}_*(\ol D\cap M)$.
This ends the proof of this case by definition of the operator
$R_{M}$.

Assume now that $p=1$. Let $f\in\ac^{1+\alpha}_*(\ol D\cap
M)$, i.e. $f\in\cc^{(1+\alpha)/2}_*(\ol D\cap M)$ and $X_\cb
f\in\ac^{\alpha}_*(\ol D\cap M)$, where $X_\cb$ is complex tangent to
$M$.

If we proceed as in the case $p=0$, but using Lemma
\ref{estimbeta-gradient} at the place of Lemma
\ref{estimbeta}, we can prove that, for each $I\in\ic'(k)$ and each
vector field $X$ tangent to $M$,
$\nabla_\zeta \wt C_{I*}(\chi X^z f)\in\bc^{1-\alpha/2}_*(D_I^*)$,
which implies by  the classical
Hardy-Littlewood lemma that $\wt C_{I*}(\chi X^z
f)\in\cc^{\alpha/2}_*(\ol D\cap M)$.

Using the case $p=0$, we get that if $X_\cb f\in\ac^{\alpha}_*(\ol
D\cap M)$, then $\wt R_{M} \chi X_\cb f\in\ac^{1+\alpha}_*(\ol D\cap
M)$. Moreover it follows from the proof of Theorem 5.6 in
\cite{BaLa} that if $X_\cb f$ is continuous then $X_\cb\wt R_{M}
(\chi f)$ and $\wt R_{M} \chi X_\cb f$ have the same regularity,
consequently $X_\cb\wt R_{M} (\chi f)\in\ac^{1+\alpha}_*(\ol D\cap
M)$.
By definition of the space $\ac^{2+\alpha}_*(\ol D\cap M)$, we have then
proved that $\wt R_{M} (\chi f)\in\ac^{2+\alpha}_*(\ol D\cap M)$.

Since by the proof
of Theorem 5.6 in \cite{BaLa}, if
$X f$ is continuous for any vector field $X$ tangent to $M$, then $X\wt R_{M}
f$ and $\wt R_{M} \chi X f$ have the
same regularity, the theorem follows for $p\geq 2$, by a simple
induction, from the definition of the spaces $\ac^{p+\alpha}_*(\ol D\cap M)$.
\end{proof}

\begin{rem}
Note that if we exchange the role played by the variables $z$ and
$\zeta$ in the operator $\wt R_{M}$, then Theorem \ref{estimlocale} is
valid for $0\leq p\leq 2l-6$
\end{rem}

Theorem \ref{estimlocale} implies better regularity properties for the
operators involve in Proposition \ref{estmeeslocales} and Theorem
\ref{formuleglobale}.

\begin{cor}\label{estimanis}
Under the hypotheses of Proposition \ref{estmeeslocales} and Theorem
\ref{formuleglobale}, the operators $T_r$, $1\leq r\leq q$ and
$n-k-q+1\le r\le n-k$, in Proposition \ref{estmeeslocales} are
bounded linear operators from $\ac^{p+\alpha}_{n,r}(M)$ into
$\ac^{p+1+\alpha}_{n,r-1}(\ol U)$ and the operators $A_r$, $1\leq
r\leq q$ and $n-k-q+1\le r\le n-k$, in  Theorem \ref{formuleglobale}
are bounded linear operators from $\ac^{p+\alpha}_{n,r}(M,E)$ into
$\ac^{p+1+\alpha}_{n,r-1}(M,E)$.
\end{cor}
\begin{proof}
Since the operators $T_r$ from  Proposition \ref{estmeeslocales} are
defined as the restrictions to $\cc^0_{n,r}(M)$ of the operators
$\wt R_{M}$, if $1\le r\le q-1$, and as the restrictions to
$\cc^0_{n,r}(M)$ of the operators $\wt R_{M}$ after exchanging the
variables $z$ and $\zeta$ in the kernels, if $n-k-q+1\le r\le n-k$,
they are in fact bounded from $\ac^{p+\alpha}_{n,r}(M)$ into
$\ac^{p+1+\alpha}_{n,r-1}(\ol U)$.

Assume $M$ is a compact $q$-concave generic CR submanifold of a
complex manifold $X$ and $E$ is an holomorphic vector bundle on X.
Let us go back to the construction of the global operators $A_r$
(cf. \cite{LaLeglobal} and \cite{BroLaLe}). First by globalizing the
homotopy formula \ref{homotopielocale} by mean of a partition of
unity, one get new linear operators $\wt T_r$ bounded from
$\ac^{p+\alpha}_{n,r}(M,E)$ into $\ac^{p+1+\alpha}_{n,r-1}(M,E)$ and
$K_r$ bounded from $\ac^{p+\alpha}_{n,r}(M,E)$ into
$\ac^{p+1+\alpha}_{n,r}(M,E)$ such that
\begin{equation}\label{homotopiecompact}
f+K_rf=\begin{cases}\wt T_{1}\opa_b f&\text{if }r=0\,,\\\opa_b\wt
T_rf+\wt T_{r+1}\opa_b f\qquad&\text{if }1\leq r\leq r-1 \text{ or }
n-k-q+1\le r\le n-k\,.
\end{cases}
\end{equation}
From functional analysis arguments we get linear operators $T'_r$,
$1\leq r\leq q$ and $n-k-q+1\le r\le n-k$, bounded from
$\cc^0_{n,r}(M,E)$ into $\cc^\infty_{n,r-1}(M,E)$ and $K'_r$, $0\leq
r\leq q-1$ and $n-k-q+1\le r\le n-k$, bounded from
$\cc^0_{n,r}(M,E)$ into $\cc^\infty_{n,r}(M,E)$ such that if
\begin{equation}
N_r=I+K_r+K'_r
\end{equation}
then
$$\cc^0_{n,r}(M,E)={\rm Im}N_r\oplus{\rm Ker}N_r,~{\rm Ker}N_r\subset
\cc^\infty_{n,r}(M,E).$$ We set $\widehat N_r=N_r+P_r=I+R_r$, where
$P_r$ denotes the linear projection with ${\rm Im}P_r={\rm Ker}N_r$
and ${\rm Ker}P_r={\rm
  Im}N_r$.
The operator $R_r$ is continuous from $\ac^{p+\alpha}_{n,r}(M,E)$
into $\ac^{p+1+\alpha}_{n,r}(M,E)$ and hence compact as an
endomorphism of $\ac^{p+\alpha}_{n,r}(M,E)$. Consequently $\widehat
{N_r}_{|_{\ac^{p+\alpha}_{n,r}(M,E)}}$ is a Fredholm endomorphism
with index $0$ and as ${\rm Ker}\widehat N_r=\{0\}$, it is an
isomorphism. Moreover $\widehat N_r^{-1}= I-R_r\widehat N_r^{-1}$,
this implies that $\widehat
{N_r^{-1}}_{|_{\ac^{p+1+\alpha}_{n,r}(M,E)}}$ is a continuous
endomorphism of $\ac^{p+1+\alpha}_{n,r}(M,E)$. The operators $A_r$
are then defined by
\begin{equation}
A_r=\begin{cases} &\widehat N_{r-1}^{-1}(T_r+T_r')\,,\qquad 1\le
r\le q \\
&(T_r+T_r')\widehat N_{r}^{-1}\,,\qquad n-k-q+1\le r\le n-k\,,
\end{cases}
\end{equation}
and hence continuous from $\ac^{p+\alpha}_{n,r}(M,E)$ into
$\ac^{p+1+\alpha}_{n,r}(M,E)$.
\end{proof}

%% file: an2.tex
\section{Second anisotropic estimates}\label{s3}

When they studied the tangential Cauchy-Riemann complex on the
Heisenberg group and more generally on strictly pseudoconvex CR
manifolds, Folland and Stein have introduced some anisotropic H\"older
spaces. 

Let $M$ be a generic $CR$ manifold of class $\cc^3$ and of real
codimension $k$ in a complex manifold $X$ of complex dimension $n$
and $D$ be a relatively compact domain in $X$.  Let
$X_1,\dots,X_{2n-2k}$ be a real basis of 
$HM$. A $\cc^1$ curve
$\gamma~:~[0,r]\to M$ is called admissible if for every $t\in [0,r]$,
\begin{equation*}
\frac{d\gamma}{dt}(t)=\sum_{j=1}^{2n-2k}c_j(t)X_j(\gamma(t))
\end{equation*}
where $\sum|c_j(t)|^2\leq 1$.

The Folland-Stein anisotropic  
H\"older spaces 
$\Gamma^{p+\alpha}(\ol D\cap M)$  are defined in the
following way:

- $\Gamma^\alpha(\ol D\cap M)$, $0<\alpha<1$, is the set of
continuous fonctions in $\ol D\cap M$ such that if for every
$x_0\in\ol D\cap M$
$$\sup_{\gamma(.)}\frac{|f(\gamma(t)-f(x_0)|}{|t|^\alpha}<\infty$$
for any admissible complex tangent curve $\gamma$ through $x_0$.

- $\Gamma^{p+\alpha}(\ol D\cap M)$, $p\ge 1$, $0<\alpha<1$, is the
set of continuous fonctions in $M$ such that $X_\cb
f\in\Gamma^{p-1+\alpha}(\ol D\cap M)$, for all complex tangent
vector fields $X_\cb$ to $M$.

If $M$ is $q$-concave, $q\ge 1$, the complex tangent vector fields and
their first Lie brackets generates all the vector fields in a
neighborhood of each point. Associated with the vector fields
$X_1,\dots,X_{2n-2k}$ we define a distance function 
${\rm dist}(x,y)$ for any $x,y\in M$ to be the
infimum of the set of all $r$ for which there exists an admissible
curve with $\gamma(0)=x$ and $\gamma(r)=y$. This distance function defines a
family of nonisotropic balls $B_r(x_0)=\{x\in M~|~{\rm
    dist}(x_0,x)<r\}$ for each point $x_0\in M$. Let $T_1,\dots,T_k$
be $k$ real vector fields such that
$X_1,\dots,X_{2n-2k},T_1,\dots,T_k$ span the whole tangent space in a
neighborhood of $x_0$. Then for $r$ sufficiently small, the ball
$B_r(x_0)$ has length comparable to $r$ in the direction of
$X_1,\dots,X_{2n-2k}$ and length comparable to $r^2$ in the direction of
$T_1,\dots,T_k$ (cf. \cite{NaStWa}). 

A function $u$ belongs to $\Gamma^\alpha(\ol D\cap M)$ if $u$ is
continuous on $\ol D\cap M$ and
\begin{equation*}
|u(x)-u(y)|\leq C ({\rm dist}(x,y))^\alpha\quad {\rm for~every}~x,y\in
\ol D\cap M 
\end{equation*}

%Note that if $M$ is $q$-concave, $q\ge 1$, the spaces
%$\Gamma^{p+\alpha}(\mb)$ are subspaces of the spaces
%$\ac^{p+\alpha}(\ol D\cap M)$ of the previous section
%for all $p\ge 0$ since for any compact set $K$ of $M$ included in a
%chart domain of $X$ there exist some
%constant $C_1$ and $C_2$ such that if $x,y\in K$
%\begin{equation*}
%C_1|x-y|\leq  {\rm dist}(x,y)\leq C_2|x-y|^\frac{1}{2}.
%\end{equation*} 

The Folland-Stein anisotropic H\"older space of forms
$\Gamma^{p+\alpha}_*(\ol D\cap M)$, $p\geq 0$, $0<\alpha<1$, is then
the space of continuous forms on $\ol D\cap M$, whose coefficients
are in $\Gamma^{p+\alpha}(\ol D\cap M)$.

We will now prove sharp estimates for the operators $\wt C_{I*}$,
$I\in\ic'(k)$, defined in Section \ref{s2} with respect to these
spaces, for this we may assume that $M$ is embedded in $\cb^n$.

\begin{prop}\label{estimationFS}
Assume $M$ is of class $\ci$ then the operators $\wt C_{I*}$,
$I\in\ic'(k)$, are continuous from 
$\Gamma^{p+\alpha}_*(\ol D\cap M)$ into $\Gamma^{p+1+\alpha}_*(\ol
D\cap M)$, $p\geq 0$.
\end{prop}
\begin{rem}
Note that in fact to get Propositon \ref{estimationFS} for some $p$
the manifold $M$ needs only to be of class $\cc^{p+4}$ and  to be of
class $\cc^{p+2}$ if we exchange the
role of $z$ and $\zeta$ in the kernel $C_{I*}$. 
\end{rem}
\begin{proof}
For $f\in\cc_{n,r}(\ol D)$, $0\leq r\leq n-k$, we set
\begin{equation*}
F(\zeta)=\wt C_{I*} f(\zeta)=\int_{z\in D\cap M}f(z)\wedge
C_{I*}(z,\zeta)=\int_{(z,\lambda)\in D\cap
  M\times\Delta_{I*}}f(z)\wedge  K_{I*}(z,\zeta,\lambda).
\end{equation*}

Without loss of generality we can assume that $f=\wt f\theta$ with
$\wt f\in \cc_{0,0}(\ol D)$ and $\theta\in\ci_{n,r}(\ol D)$. We set 
$\psi(\zeta)=\wt C_{I*} \theta(\zeta)$, the differential form
$\psi$ is then of class $\ci$ (cf. Section 5 in \cite{BaLa}) and for $\zeta_1$ and $\zeta_2$ in $\ol D$
\begin{equation*}
F(\zeta_1)-F(\zeta_2)=\int_{z\in D\cap M}(\wt f(z)-\wt f(\zeta_1))\theta(z)\wedge
(C_{I*}(z,\zeta_1)-C_{I*}(z,\zeta_2)+\wt f(\zeta_1)(\psi(\zeta_1)-\psi(\zeta_2)).
\end{equation*}

We first concider the case where $p=0$.
Let $\gamma$ be an admissible curve with $\gamma(0)=\zeta_1$ and
$\gamma(r)=\zeta_2$ we have to estimate $\frac{d(F\circ\gamma)}{dt}(0)-
\frac{d(F\circ\gamma)}{dt}(r)$ when $\wt f\in\Gamma^{\alpha}(\ol D\cap
M)$. We have
\begin{align*}
\frac{d(F\circ\gamma)}{dt}(0)&-
\frac{d(F\circ\gamma)}{dt}(r)\\
&=\int_{z\in D\cap M}(\wt f(z)-\wt
f(\gamma(0)))\theta(z)\wedge (\frac{dC_{I*}(z,\gamma(.))}{dt}(0)- 
\frac{dC_{I*}(z,\gamma(.))}{dt}(r))\\
&+\wt
f(\gamma(0))(\frac{d\psi\circ\gamma}{dt}(0)- \frac{d\psi\circ\gamma}{dt}(r)).
\end{align*}

We set $D_1=\{z\in\ol D\cap M~|~|z-\zeta_1|\leq 2|\zeta_1-\zeta_2|\}$
and $D_2=\{z\in\ol D\cap M~|~|z-\zeta_1|\geq 2|\zeta_1-\zeta_2|\}$
then $\ol D\cap M=D_1\cup D_2$. We then have 
\begin{equation*}
|\frac{d(F\circ\gamma)}{dt}(0)-
\frac{d(F\circ\gamma)}{dt}(r)|\leq J_1+J_2+|\wt
f(\gamma(0))(\frac{d\psi\circ\gamma}{dt}(0)-
\frac{d\psi\circ\gamma}{dt}(r))|
\end{equation*}
where
\begin{equation}\label{estimJ1}
\begin{aligned}
J_1&=|\int_{z\in D_1}(\wt f(z)-\wt
f(\gamma(0)))\theta(z)\wedge \big(\frac{dC_{I*}(z,\gamma(.))}{dt}(0)- 
\frac{dC_{I*}(z,\gamma(.))}{dt}(r)\big)|\\
&\leq\int_{z\in D_1}|\wt f(z)-\wt
f(\gamma(0))||\theta(z)|\wedge |\frac{dC_{I*}(z,\gamma(.))}{dt}(0)|\\
&+\int_{z\in D_1}|\wt f(z)-\wt
f(\gamma(r))||\theta(z)|\wedge |\frac{dC_{I*}(z,\gamma(.))}{dt}(r)|\\
&+|\wt f(\gamma(0))-\wt f(\gamma(r))|\int_{z\in D_1}|\theta(z)|\wedge
|\frac{dC_{I*}(z,\gamma(.))}{dt}(r)|. 
\end{aligned}
\end{equation}
and
\begin{align*}
J_2&=|\int_{z\in D_2}(\wt f(z)-\wt
f(\gamma(0)))\theta(z)\wedge (\frac{dC_{I*}(z,\gamma(.))}{dt}(0)- 
\frac{dC_{I*}(z,\gamma(.))}{dt}(r))|\\
&\leq \int_{z\in D_2}|\wt f(z)-\wt
f(\gamma(0))||\theta(z)|\wedge |\frac{dC_{I*}(z,\gamma(.))}{dt}(0)- 
\frac{dC_{I*}(z,\gamma(.))}{dt}(r)|.
\end{align*}

Let us first consider $J_1$. The last term in (\ref{estimJ1}) is
clearly controlled by $C\|\wt f\|_{\Gamma^\alpha}r^\alpha$ and both other integral
terms in (\ref{estimJ1}) are of the same type.

Note that since $\gamma$ is a complex tangent curve
$\frac{dC_{I*}(z,\gamma(.))}{dt}(s)$, $s=0,r$, involves only complex tangential
derivatives of the part of degree $k$ in $\lambda$ of the kernel
$K_{I*}(z,\zeta,\lambda)$. As $X_\cb^\zeta \phi(z,\zeta,\lambda)=O_1$
the coefficients of these derivatives are sums of terms of type
$\frac{O_{k-m}}{\phi^n}$, $0\leq m\leq k$. Moreover the first order
terms in $\Phi(z,\zeta,\lambda)$ are equivalent to the distance of $z$
and $\zeta$ in the real tangential directions transverse to the complex
tangent space (i.e. in the directions span by $T_1,\dots,T_k$ for a
good choice of $T_1,\dots,T_k$). Therefore $[{\rm dist}(z,\zeta)]^2=
O(|\Phi(z,\zeta,\lambda)|)$ and if $\wt f\in\Gamma^\alpha(\ol
D\cap M)$ we get for $s=0,r$
\begin{equation*}
|\wt f(z)-\wt f(\gamma(s))|=O(|\Phi(z,\gamma(s),\lambda)|^{\alpha/2}).
\end{equation*}

The same calculations as in Section 5.1 of \cite{BaLa} lead to the
estimate 
$$J_1\leq C\|\wt f\|_{\Gamma^\alpha}r^\alpha.$$

To estimate $J_2$, it follows from Section 5.1 of \cite{BaLa} that the
main point is to control the difference
$$\frac{1}{\Phi(z,\gamma(0),\lambda)}-\frac{1}{\Phi(z,\gamma(r),\lambda)}.$$
Since $[d_\zeta\Phi(z,\gamma(t),\lambda).\gamma'(t)]_{|t=r}=O_1$ we
have 
\begin{equation*}
|\frac{1}{\Phi(z,\gamma(0),\lambda)}-
\frac{1}{\Phi(z,\gamma(r),\lambda)}|\leq C    
\sum_{p=0}^{n-1} \frac{|\gamma(0)-\gamma(r)|O_1+
  |\gamma(0)-\gamma(r)|^2O_0}{
  |\Phi^{n-p}(z,\gamma(0),\lambda)||\Phi^{p+1}(z,\gamma(r),\lambda)|}.
\end{equation*}
As for $J_1$ following the estimations in Section 5.1 of \cite{BaLa} we
get 
$$J_2\leq C\|\wt f\|_{\Gamma^\alpha}r^\alpha$$
and finally
$$|\frac{d(F\circ\gamma)}{dt}(0)-
\frac{d(F\circ\gamma)}{dt}(r)|\leq C\|\wt f\|_{\Gamma^\alpha}r^\alpha$$
since $\wt f$ is bounded and $\psi$ is of class $\cc^1$.

Using Lemma \ref{derivation} we can derive the case $p\geq 1$ from the
case $p=0$ in the same way as in the previous section and in Section
5.2 in \cite{BaLa}.
\end{proof}

It follows from Proposition \ref{estimationFS} and from the definition of
 the  operators $T_r$
and $A_r$, $1\leq r\leq q$ and $n-k-q+1\leq r\leq n-k$, from Section
\ref{s1} that Corollary 
\ref{estimanis} still holds if we replace the spaces
$\ac^{p+\alpha}_{n,r}(M)$ and $\ac^{p+1+\alpha}_{n,r-1}(\ol U)$ by the
spaces  $\Gamma^{p+\alpha}_{n,r}(M)$ and
$\Gamma^{p+1+\alpha}_{n,r-1}(\ol U)$ respectively and the spaces
$\ac^{p+\alpha}_{n,r}(M,E)$ and $\ac^{p+1+\alpha}_{n,r-1}(M,E)$ by the
spaces  $\Gamma^{p+\alpha}_{n,r}(M,E)$ and
$\Gamma^{p+1+\alpha}_{n,r-1}(M,E)$
respectively.